\documentclass[runningheads]{llncs}
\usepackage{graphicx}
 \usepackage{color}
 \usepackage{amssymb, latexsym, amsmath, epsfig}
\begin{document}
\title{Outlier removal for isogeometric spectral approximation with the optimally-blended quadratures \thanks{This work and visit of Quanling Deng in Krakow was partially supported by National Science Centre, Poland grant no. 017/26/M/ST1/00281. This publication was made possible in part by the CSIRO Professorial Chair in Computational Geoscience at Curtin University and the Deep Earth Imaging Enterprise Future Science Platforms of the Commonwealth Scientific Industrial Research Organisation, CSIRO, of Australia. Additional support was provided by the European Union's Horizon 2020 Research and Innovation Program of the Marie Sklodowska-Curie grant agreement No. 777778 and the Curtin Institute for Computation.}}
\titlerunning{Outlier removal in IGA with blended quadratures}
%
\author{Quanling Deng\inst{1}\orcidID{0000-0002-6159-1233} \and
Victor M. Calo\inst{2}\orcidID{0000-0002-1805-4045}}
\authorrunning{Q. Deng and V. Calo.}
%
\institute{Department of Mathematics, University of Wisconsin--Madison, Madison, WI 53706, USA \email{Quanling.Deng@math.wisc.edu} \and
Curtin Institute for Computation \& School of Electrical
    Engineering, Computing and Mathematical Sciences, Curtin
    University, P.O. Box U1987, Perth, WA 6845, Australia \email{Victor.Calo@curtin.edu.au}}
\maketitle              
\begin{abstract}
It is well-known that outliers appear in the high-frequency region in the approximate spectrum of 
isogeometric analysis of the second-order elliptic operator. 
Recently, the outliers have been eliminated by a boundary penalty technique.
The essential idea is to impose extra conditions arising from the differential equation at the domain boundary.  
In this paper, we extend the idea to remove outliers in the superconvergent approximate spectrum of  
isogeometric analysis with optimally-blended quadrature rules. 
We show numerically that the eigenvalue errors are of superconvergence rate $h^{2p+2}$ 
and the overall spectrum is outlier-free. 
The condition number and stiffness of the resulting algebraic system are reduced significantly. Various numerical examples demonstrate the performance of the proposed method.

\end{abstract}
\section{Introduction} \label{sec:intr} 

Isogeometric analysis (IGA) is a widely-used analysis tool that combines the classical finite element analysis with computer-aided design and analysis tools. It was introduced in 2005~\cite{hughes2005isogeometric,cottrell2009isogeometric}. There is a rich literature since its first development; see an overview paper \cite{nguyen2015isogeometric} and the references therein. 
In particular, a rich literature on IGA has been shown that the method outperforms the 
classical finite element method (FEM) on the spectral approximations of the second-order elliptic operators. 
The initial work \cite{cottrell2006isogeometric} showed that the spectral errors of IGA were significantly smaller when compared with FEM approximations.
In~\cite{hughes2014finite}, the authors further explored the advantages of IGA on the spectral approximations.

To further reduce the spectral errors on IGA spectral approximations, 
on one hand, 
the recent work \cite{calo2019dispersion,puzyrev2017dispersion} introduced Gauss-Legendre and 
Gauss-Lobatto optimally-blended quadrature rules. 
By invoking the dispersion analysis which was unified with the spectral analysis in \cite{hughes2008duality}, 
the spectral errors were shown to be superconvergent with two extra orders. 
The work \cite{deng2018ddm} generalized the blended rules to arbitrary $p$-th order IGA with maximal continuity. 
Along the line, the work \cite{deng2018dispersion} further studied the computational efficiency and the work
\cite{bartovn2018generalization,calo2017quadrature,puzyrev2018spectral,deng2019optimal,deng2018isogeometric} studied its applications.

On the other hand, the spectral errors in the highest-frequency regions are much larger than in the lower-frequency region. 
There is a thin layer in the highest-frequency region which is referred to as ``outliers". 
The outliers in isogeometric spectral approximations were first observed in \cite{cottrell2006isogeometric} in 2006.  
The question of how to efficiently remove the outliers remained open until recently. 
In \cite{deng2020boundary}, the authors  removed the outliers by a boundary penalty technique. 
The main idea is to impose higher-order consistency conditions on the boundaries to the isogeometric spectral approximations. 
The work proposed to impose these conditions weakly. 
Outliers are eliminated and the condition numbers of the systems are reduced significantly. 

In this paper, we propose to further reduce spectral errors by combining the optimally-blended quadrature rules 
and the boundary penalty technique. 
To illustrate the idea, we focus on tensor-product meshes on rectangular domains. 
We first develop the method in 1D and obtain the generalized matrix eigenvalue problems. 
We then apply the tensor-product structure to generate the matrix problems in multiple dimensions.
By using the optimally-blended quadrature rules, we retain the eigenvalue superconvergence rate.
By applying the boundary penalty technique, we remove the outliers in the superconvergent spectrum. 
The method also reduces the condition numbers.

The rest of this paper is organized as follows. 
Section \ref{sec:ps} presents the problem and its discretization by the standard isogeometric analysis. 
Section \ref{sec:ob} concerns the Gauss-Legendre and Gauss-Lobatto quadrature rules. We then present the 
optimally-blended rules. 
In Section \ref{sec:mi}, we apply the boundary penalty technique developed in \cite{deng2020boundary} to the IGA setting with blending rules. 
Section \ref{sec:num} collects numerical results that demonstrate the performance of the proposed method. 
In particular, we perform the numerical study on the condition numbers.
Concluding remarks are presented in Section \ref{sec:conclusion}.

\section{Problem setting} \label{sec:ps}
Let $\Omega = [0,1]^d \subset  \mathbb{R}^d, d=1,2,3$ be a bounded open domain 
with Lipschitz boundary $\partial \Omega$. 
We use the standard notation for the Hilbert and Sobolev spaces. 
For a measurable subset $S\subseteq \Omega$, 
we denote by $(\cdot,\cdot)_S$ and $\| \cdot \|_S$ the $L^2$-inner product and its norm, respectively. 
We omit the subscripts when it is clear in the context.
For any integer $m\ge1$,
we denote the $H^m$-norm and $H^m$-seminorm as $\| \cdot \|_{H^m(S)}$ and 
$| \cdot |_{H^m(S)}$, respectively. In particular, we denote by $H^1_0(\Omega)$ the Sobolev space 
with functions in $H^1(\Omega)$ that are vanishing at the boundary.
We consider the classical second-order elliptic eigenvalue problem: 
Find the eigenpairs $(\lambda, u)\in \mathbb{R}^+\times H^1_0(\Omega)$ with $\|u\|_\Omega=1$ such that
\begin{equation} \label{eq:pde}
\begin{aligned}
- \Delta u & =  \lambda u \quad &&  \text{in} \quad \Omega, \\
u & = 0 \quad && \text{on} \quad  \partial \Omega,
\end{aligned}
\end{equation}
where $\Delta = \nabla^2$ is the Laplacian. 
The variational formulation of \eqref{eq:pde} is to find $\lambda \in \mathbb{R}^{+}$ and 
$u \in H^1_0(\Omega)$ with $\|u\|_\Omega=1$ such that 
\begin{equation} \label{eq:vf}
a(w, u) =  \lambda b(w, u), \quad \forall \ w \in H^1_0(\Omega), 
\end{equation}
where the bilinear forms
\begin{equation}
a(v,w) := (\nabla v, \nabla w)_\Omega,
\qquad
b(v,w) := (v,w)_\Omega.
\end{equation}

It is well-known that the eigenvalue problem \eqref{eq:vf} has a countable set of positive eigenvalues 
(see, for example, \cite[Sec. 9.8]{Brezis:11})
\begin{equation*}
0 < \lambda_1 < \lambda_2 \leq \lambda_3 \leq \cdots
\end{equation*}
and an associated set of orthonormal eigenfunctions $\{ u_j\}_{j=1}^\infty$, that is, 
$
(u_j, u_k) = \delta_{jk}, 
$
where $\delta_{jk} =1$ is the Kronecker delta. 
Consequently, the eigenfunctions are also orthogonal in the energy inner product since there holds 
$
a(u_j, u_k) =  \lambda_j b(u_j, u_k) = \lambda_j \delta_{jk}.
$

At the discretize level, we first discretize the domain $\Omega$ with a uniform tensor-product mesh. 
We denote a general element as $\tau$ and its collection as $\mathcal{T}_h$ such that 
$\overline \Omega = \cup_{\tau\in \mathcal{T}_h}  \tau$. 
Let $h = \max_{\tau \in \mathcal{T}_h} \text{diameter}(\tau)$. 
In the IGA setting, for simplicity, we use the B-splines. 
The B-spline basis functions in 1D are defined by using the Cox-de Boor recursion formula; 
we refer to \cite{de1978practical,piegl2012nurbs} for details.
Let $X = \{x_0, x_1, \cdots, x_m \}$ be a knot vector with knots $x_j$, that is, 
a nondecreasing sequence of real numbers.  
The $j$-th B-spline basis function of degree $p$, denoted as $\phi^j_p(x)$, is defined recursively as 
\begin{equation} \label{eq:Bspline}
\begin{aligned}
\phi^j_0(x) & = 
\begin{cases}
1, \quad \text{if} \ x_j \le x < x_{j+1}, \\
0, \quad \text{otherwise}, \\
\end{cases} \\ 
\phi^j_p(x) & = \frac{x - x_j}{x_{j+p} - x_j} \phi^j_{p-1}(x) + \frac{x_{j+p+1} - x}{x_{j+p+1} - x_{j+1}} \phi^{j+1}_{p-1}(x).
\end{aligned}
\end{equation}
A tensor-product of these 1D B-splines produces the B-spline basis functions in multiple dimensions. 
We define the multi-dimensional approximation space as $V^h_p \subset H^1_0(\Omega)$ with (see \cite{buffa2010isogeometric,evans2013isogeometric} for details):
\begin{equation*}
V^h_p = \text{span} \{ \phi_j^p \}_{j=1}^{N_h} = 
\begin{cases}
 \text{span} \{ \phi^{j_x}_{p_x}(x) \}_{j_x=1}^{N_x}, & \text{in 1D}, \\
 \text{span} \{ \phi^{j_x}_{p_x}(x) \phi^{j_y}_{p_y}(y) \}_{j_x, j_y=1}^{N_x, N_y}, & \text{in 2D}, \\
 \text{span} \{\phi^{j_x}_{p_x}(x) \phi^{j_y}_{p_y}(y) \phi^{j_z}_{p_z}(z) \}_{j_x, j_y, j_z=1}^{N_x, N_y,N_z}, & \text{in 3D}, \\
\end{cases}
\end{equation*}
where $p_x, p_y, p_z$ specify the approximation order in each dimension. 
$N_x, N_y, N_z$ is the total number of basis functions in each dimension 
and $N_h$ is the total number of degrees of freedom.
The isogeometric analysis of \eqref{eq:pde} in variational formulation seeks $\lambda^h \in \mathbb{R}$ and $u^h \in V^h_p$ with $\| u^h \|_\Omega = 1$ such that 
\begin{equation} \label{eq:vfh}
a(w^h, u^h) =  \lambda^h b(w^h, u^h), \quad \forall \ w^h \in V^h_p.
\end{equation}

At the algebraic level, we approximate the eigenfunctions as a linear combination of the B-spline basis functions 
and substitute all the B-spline basis functions for $w^h$ in \eqref{eq:vfh}. 
This leads to the generalized matrix eigenvalue problem
\begin{equation} \label{eq:mevp}
\mathbf{K} \mathbf{U} = \lambda^h \mathbf{M} \mathbf{U},
\end{equation}
where $\mathbf{K}_{kl} =  a(\phi_p^k, \phi_p^l), \mathbf{M}_{kl} = b(\phi_p^k, \phi_p^l),$ 
and $\mathbf{U}$ is the corresponding representation of the eigenvector as the coefficients of the B-spline basis functions. 
In practice, we evaluate the integrals involved in the bilinear forms $a(\cdot, \cdot) $ and $b(\cdot, \cdot)$ 
numerically using quadrature rules. 
In the next section, we present the Gauss--Legendre and Gauss--Lobatto quadrature rules, 
followed by their optimally-blended rules.

\section{Quadrature rules and optimal blending} \label{sec:ob}
In this section, we first present the classic Gauss-type quadrature rules and 
then present the optimally-blended rules developed recently in \cite{calo2019dispersion}. 
While the optimally-blended rules have been developed in \cite{deng2018ddm} for arbitrary 
order isogeometric elements, we focus on the lower-order cases for simplicity. 

\subsection{Gaussian quadrature rules}
Gaussian quadrature rules are well-known and 
we present these rules by following the book \cite{kythe2004handbook}.
On a reference interval $\hat \tau$, a quadrature rule is of the form
\begin{equation} \label{eq:qr}
\int_{\hat \tau} \hat f(\hat{\boldsymbol{x}}) \ \text{d} \hat{\boldsymbol{x}} \approx \sum_{l=1}^{m} \hat{\varpi}_l \hat f (\hat{n_l}),
\end{equation}
where $\hat{\varpi}_l$ are the weights, $\hat{n_l}$ are the nodes, 
and $m$ is the number of quadrature points. 
We list the lower-order Gauss-Legendre and Gauss-Lobatto quadrature rules in 1D below; 
we refer to \cite{kythe2004handbook} for rules with more points. 
The Gauss-Legendre quadrature rules for $m=1,2,3,4$  in the reference interval $[-1, 1]$ are as follows:
\begin{equation} \label{eq:gquad}
\begin{aligned}
m=1: \qquad & \hat n_1 = 0, \qquad \hat \varpi_1  = 2; \\
m=2: \qquad & \hat n_1 = \pm \frac{\sqrt{3}}{3}, \qquad \hat \varpi_1  = 1; \\
m=3: \qquad & \hat n_1  = 0, \quad \hat n_{2,3} = \pm \sqrt{\frac{3}{5}}, \qquad \hat \varpi_1 = \frac{8}{9}, \quad \hat \varpi_{2,3}  = \frac{5}{9}; \\
m=4: \qquad & \hat n_{1,2,3,4} = \pm \sqrt{ \frac{3}{7} \mp \frac{2}{7} \sqrt{\frac{6}{5}} }, \qquad \hat \varpi_{1,2,3,4}  = \frac{18\pm\sqrt{30} }{36}. \\
\end{aligned}
\end{equation}
A Gauss-Legendre quadrature rule with $m$ points integrates exactly a polynomial of degree $2m-1$ or less.
The Gauss-Lobatto quadrature rules with $m=2,3,4,5$ in the reference interval $[-1, 1]$ are as follows: 
\begin{equation} \label{eq:lquad}
\begin{aligned}
m=2: \qquad & \hat n_{1,2} = \pm 1, \qquad \hat \varpi_1  = 1; \\
m=3: \qquad & \hat n_1 = 0, \quad \hat n_{2,3} =  \pm 1, \qquad \hat \varpi_1  = \frac{4}{3}, \quad \hat \varpi_{2,3}  = \frac{1}{3}; \\
m=4: \qquad & \hat n_{1,2} = \pm \sqrt{\frac{1}{5}}, \quad \hat n_{3,4}  = \pm 1, \quad  \qquad \hat \varpi_{1,2} = \frac{5}{6}, \quad \hat \varpi_{3,4}  = \frac{1}{6}; \\
m=5: \qquad & \hat n_1 = 0, \quad \hat n_{2,3} =  \pm \sqrt{\frac{3}{7}},  \quad \hat n_{4,5} = \pm 1,  \quad \hat \varpi_1  = \frac{32}{45}, \quad \hat \varpi_{2,3} = \frac{49}{90}, \hat \varpi_{4,5}  = \frac{1}{10}. 
\end{aligned}
\end{equation}
A Gauss-Lobatto quadrature rule with $m$ points integrates exactly a polynomial of degree $2m-3$ or less.

For each element $\tau$, there is a one-to-one and onto mapping $\sigma$ such that 
$\tau = \sigma(\hat \tau)$, which leads to the correspondence between the functions 
on $\tau$ and $\hat \tau$. Let $J_\tau$ be the corresponding Jacobian of the mapping. Using the mapping, 
\eqref{eq:qr} induces a quadrature rule over the element $\tau$ given by
\begin{equation} \label{eq:q}
\int_{\tau}  f(\boldsymbol{x}) \ \text{d} \boldsymbol{x} \approx \sum_{l=1}^{N_q} \varpi_{l,\tau} f (n_{l,\tau}),
\end{equation}
where $\varpi_{l,\tau} = \text{det}(J_\tau) \hat \varpi_l$ and $n_{l,\tau} = \sigma(\hat n_l)$. 
For blending rules, we denote by $Q$ a general quadrature rule, 
by $G_m$ the $m-$point Gauss-Legendre quadrature rule, 
by $L_m$ the $m-$point Gauss-Lobatto quadrature rule, 
by $Q_\eta$ a blended rule, and 
by $O_p$ the optimally-blended rule for the $p$-th order isogeometric analysis.

\subsection{Optimal blending quadrature rules}
Let $Q_1 = \{\varpi_{l,\tau}^{(1)}, n_{l,\tau}^{(1)} \}_{l=1}^{m_1}$ and 
$Q_2 = \{\varpi_{l,\tau}^{(2)}, n_{l,\tau}^{(2)} \}_{l=1}^{m_2}$ be two quadrature rules. 
We define the blended quadrature rule as
\begin{equation}
Q_\eta = \eta Q_1 + (1-\eta) Q_2,
\end{equation}
where $\eta \in \mathbb{R}$ is a blending parameter. 
We note that the blending parameter can be both positive and negative. 
The blended rule $Q_\eta$ for the integration of a function $f$ is understood as
\begin{equation}
\begin{aligned}
\int_{\tau}  f(\boldsymbol{x}) \ \text{d} \boldsymbol{x} & \approx \tau \sum_{l=1}^{m_1} \varpi_{l,\tau}^{(1)} f (n_{l,\tau}^{(1)}) + (1 - \tau) \sum_{l=1}^{m_2} \varpi_{l,\tau}^{(2)} f (n_{l,\tau}^{(2)}). \\
\end{aligned}
\end{equation}

The blended rules for isogeometric analysis of the eigenvalue problem can reduce significantly the spectral errors. 
In particular, the optimally-blended rules deliver two extra orders of convergence for the eigenvalue errors; 
see \cite{calo2019dispersion,deng2018ddm} for the developments. 
We present the optimally-blended rule for the $p$-th order isogeometric elements 
\begin{equation} \label{eq:br}
O_p = Q_\eta = \eta G_{p+1} + (1-\eta) L_{p+1},
\end{equation}
where the optimal blending parameters are given in the Table \ref{tab:ob} below.

\begin{table}[ht]
\label{tab:ob} 
\caption{Optimal blending parameters for isogeometric elements.} 
\centering 
{\renewcommand{\arraystretch}{1.5}
\begin{tabular}{| c | c | c | c | c | c | c| c| }
\hline
$p$ & 1 & 2   & 3 & 4 & 5 & 6 & 7 \\[0.0cm] \hline

$\eta$  & $\frac{1}{2}$ & $\frac{1}{3}$ & $-\frac{3}{2}$ & $-\frac{79}{5}$ & $-174$ & $-\frac{91177}{35}$ &  $-\frac{105013}{2}$ \\[0.1cm] \hline
\end{tabular}
}
\end{table}
\vspace{-0.5cm}
We remark that there are other optimally-blended quadrature rules for isogeometric elements 
developed in \cite{calo2019dispersion,deng2018ddm}. 
Non-standard quadrature rules are developed in \cite{deng2018dispersion} and 
they are shown to be equivalent with the optimally-blended rules. 
They all lead to the same stiffness and mass matrices. 
Herein, for simplicity, we only adopt the blended rules in the form of \eqref{eq:br}.

\section{The boundary penalty technique} \label{sec:mi}
Outliers appear in the isogeometric spectral approximations when using $C^2$ cubic B-splines  
and higher-order  approximations. 
These outliers can be removed by a boundary penalty technique. 
In this section, we first recall the boundary penalty technique introduced recently in \cite{deng2020boundary}. 
We present the idea for 1D problem with $\Omega = [0,1] $ and $\partial \Omega = \{0, 1\}$. 
We then generalize it at the algebraic level using the tensor-product structure for the problem in multiple dimensions. 
We denote 
\begin{equation}
\alpha = \lfloor \frac{p-1}{2} \rfloor =
\begin{cases}
\frac{p-1}{2}, & p \quad \text{is odd}, \\
\frac{p-2}{2}, & p \quad \text{is even}. \\
\end{cases}
\end{equation}
The isogeometric analysis of \eqref{eq:pde} in 1D with the boundary penalty technique is: 
Find $\tilde \lambda^h \in \mathbb{R}$ and $\tilde u^h \in V^h_p$ such that 
\begin{equation} \label{eq:dcvfh}
\tilde a(w^h, \tilde u^h) =  \tilde\lambda^h \tilde b(w^h, \tilde u^h), \quad \forall \ w^h \in V^h_p.
\end{equation} 
Herein, for $w, v \in V^h_p$
\begin{subequations} \label{eq:dcvfbfs}
\begin{align}
\tilde a(w,  v) & = \int_0^1 w' v' \ d x + \sum_{\ell = 1}^\alpha \eta_{a,\ell} \pi^2 h^{6\ell-3} \Big( w^{(2\ell)}(0) v^{(2\ell)}(0) + w^{(2\ell)}(1) v^{(2\ell)}(1) \Big), \label{eq:dcvfbfsa} \\
\tilde b(w,  v) & = \int_0^1 w v \ d x + \sum_{\ell = 1}^\alpha \eta_{b,\ell} h^{6\ell-1} \Big( w^{(2\ell)}(0) v^{(2\ell)}(0) + w^{(2\ell)}(1) v^{(2\ell)}(1) \Big), \label{eq:dcvfbfsb}
\end{align}
\end{subequations}
where $\eta_{a,\ell}, \eta_{b,\ell}$ are penalty parameters set to $\eta_{a,\ell} = \eta_{b,\ell} = 1$ in default. 
The superscript $(2\ell)$ denotes the $2\ell$-th derivative. 
We further approximate the inner-products by the quadrature rules discussed above. 
For $p$-th order element and $\tau \in \mathcal{T}_h$, 
we denote $G_{p+1} = \{\varpi_{l,\tau}^G, n_{l,\tau}^G \}_{l=1}^{p+1}$ and $L_{p+1} = \{\varpi_{l,\tau}^L, n_{l,\tau}^L \}_{l=1}^{p+1}$.
Applying the optimally-blended quadrature rules in the form of \eqref{eq:br} to \eqref{eq:dcvfh}, 
we obtain the approximated form
\begin{equation} \label{eq:vfho}
\tilde a_h(w^h, \tilde u^h) =  \tilde\lambda^h  \tilde b_h(w^h, \tilde u^h), \quad \forall \ w^h \in V^h_p,
\end{equation}
where for $w,v \in V^h_p$
\begin{equation} \label{eq:ba}
\begin{aligned}
\tilde a_h(w, v) & = \sum_{\tau \in \mathcal{T}_h} \sum_{l=1}^{p+1} \big( \eta \varpi_{l,\tau}^G \nabla w (n_{l,\tau}^G ) \cdot \nabla v (n_{l,\tau}^G ) + (1-\eta) \varpi_{l,\tau}^L \nabla w (n_{l,\tau}^L ) \cdot \nabla v (n_{l,\tau}^L ) \big) \\
& \quad + \sum_{\ell = 1}^\alpha \eta_{a,\ell} \pi^2 h^{6\ell-3} \Big(w^{(2\ell)}(0) v^{(2\ell)}(0) + w^{(2\ell)}(1) v^{(2\ell)}(1) \Big),
\end{aligned}
\end{equation}
and
\begin{equation} \label{eq:bb}
\begin{aligned}
\tilde b_h(w, v) & = \sum_{\tau \in \mathcal{T}_h} \sum_{l=1}^{p+1} \big( \eta \varpi_{l,\tau}^G w (n_{l,\tau}^G ) \cdot v (n_{l,\tau}^G ) + (1 - \eta)  \varpi_{l,\tau}^L w (n_{l,\tau}^L ) \cdot v (n_{l,\tau}^L ) \big) \\
& \quad + \sum_{\ell = 1}^\alpha \eta_{b,\ell} h^{6\ell-1} \Big(w^{(2\ell)}(0) v^{(2\ell)}(0) + w^{(2\ell)}(1) v^{(2\ell)}(1) \Big).
\end{aligned}
\end{equation}
With this in mind, we arrive at the matrix eigenvalue problem
\begin{equation} \label{eq:mevp1d}
\tilde{\mathbf{K}} \mathbf{U} = \tilde \lambda^h \tilde{\mathbf{M}} \mathbf{U},
\end{equation}
where $\tilde{\mathbf{K}}_{kl} = \tilde{\mathbf{K}}^{1D}_{kl} =  \tilde a_h(\phi_p^k, \phi_p^l), \tilde{\mathbf{M}}_{kl} = \tilde{\mathbf{M}}^{1D}_{kl} = \tilde b_h(\phi_p^k, \phi_p^l)$ in 1D. 
Using the tensor-product structure and introducing the outer-product $\otimes$ (also known as the Kronecker product), 
the corresponding 2D matrices (see \cite{calo2019dispersion} for details) are
\begin{equation} \label{eq:mevp2d}
\begin{aligned}
\tilde{\mathbf{K}} & = \tilde{\mathbf{K}}^{2D} = \tilde{\mathbf{K}}^{1D}_x \otimes \tilde{\mathbf{M}}^{1D}_y + \tilde{\mathbf{M}}^{1D}_x \otimes \tilde{\mathbf{K}}^{1D}_y, \\
\tilde{\mathbf{M}} & = \tilde{\mathbf{M}}^{2D} = \tilde{\mathbf{M}}^{1D}_x \otimes \tilde{\mathbf{M}}^{1D}_y,
\end{aligned}
\end{equation}
and the 3D matrices are
\begin{equation} \label{eq:mevp3d}
\begin{aligned}
\tilde{\mathbf{K}} & = \tilde{\mathbf{K}}^{3D} = \tilde{\mathbf{K}}^{1D}_x \otimes \tilde{\mathbf{M}}^{1D}_y \otimes \tilde{\mathbf{M}}^{1D}_z + \tilde{\mathbf{M}}^{1D}_x \otimes \tilde{\mathbf{K}}^{1D}_y \otimes \tilde{\mathbf{M}}^{1D}_z + \tilde{\mathbf{M}}^{1D}_x \otimes \tilde{\mathbf{M}}^{1D}_y \otimes \tilde{\mathbf{K}}^{1D}_z, \\
\tilde{\mathbf{M}} & = \tilde{\mathbf{M}}^{3D} = \tilde{\mathbf{M}}^{1D}_x \otimes \tilde{\mathbf{M}}^{1D}_y \otimes \tilde{\mathbf{M}}^{1D}_z,
\end{aligned}
\end{equation}
where $\tilde{\mathbf{K}}^{1D}_q,  \tilde{\mathbf{M}}^{1D}_q, q=x,y,z$ are 1D matrices generated from the modified bilinear forms in \eqref{eq:vfho}.
We remark that in the limiting case with uniform meshes, analytical eigenpairs to the matrix eigenvalue problem 
\ref{eq:mevp1d} can be derived based on \cite{deng2021}; see \cite{deng2020boundary} for IGA with Gauss quadrature rule. We omit it here for brevity.

\section{Numerical examples} \label{sec:num}
In this section, we present numerical tests to demonstrate the performance of the method.  
We consider the problem \eqref{eq:pde}  with $d=1,2,3.$
The 1D problem has true eigenpairs
$ \big(
\lambda_j = j^2 \pi^2,  u_j = \sin( j\pi x) \big), j = 1, 2, \cdots,
$
the 2D problem has true eigenpairs
$ \big(
\lambda_{jk} = ( j^2 + k^2 ) \pi^2, u_{jk} = \sin( j\pi x)\sin( k\pi y) \big),  j,k = 1, 2, \cdots,
$
and the 3D problem has true eigenpairs
$ \big(
\lambda_{jkl} = ( j^2 + k^2 + l^2) \pi^2, u_{jkl} = \sin( j\pi x)\sin( k\pi y) \sin( l\pi z) \big),  
j,k,l = 1, 2, \cdots.
$
 We sort both the exact and approximate eigenvalues in ascending order. 
Since outliers appear in the spectrum for cubic ($p=3$) and higher order isogeometric elements, 
we focus on $p$-th order elements with $p\ge3$.

\subsection{Numerical study on error convergence rates and outliers}
To study the errors, we consider both the $H^1$-seminorm and $L^2$-norm for the eigenfunctions. 
The optimal convergence rates in $H^1$-seminorm and $L^2$-norm are $h^p$ and 
$h^{p+1}$ for $p$-th order elements, respectively. For the eigenvalues, we consider the relative eigenvalue errors defined as
$
\frac{|\tilde \lambda^h_j - \lambda_j |}{\lambda_j}.
$
The optimal convergence rate is $h^{2p}$ for $p$-th order elements. We denote by $\rho_p$ the convergence rate of $p$-th order isogeometric elements.

Table \ref{tab:eig1d} shows the eigenvalue and eigenfunction errors for $p=3,4,5$ in 1D while Table \ref{tab:eig2d3d} shows the eigenvalue errors for the problems in 2D and 3D. 
The first eigenvalue error approaches the machine precision fast. 
Thus, we calculate the convergence rates with coarser meshes.  
In all these scenarios, the eigenfunction errors are convergent optimally, 
while the eigenvalue errors are superconvergent with two extra orders, i.e., $h^{2p+2}$. 
These results are in good agreement with the theoretical predictions. 

\begin{table}[ht]
\caption{Errors and convergence rates for the first and sixth eigenpairs in 1D when using IGA with optimally-blended quadratures and the boundary penalty technique.}
\label{tab:eig1d} 
\centering 
\begin{tabular}{| c | c || ccc | ccc | cc |}
\hline
$p$ & $N$ & $\frac{|\tilde \lambda^h_1-\lambda_1|}{\lambda_1}$ &  $|u_1- \tilde u^h_1|_{H^1}$ & $\| u_1- \tilde u_1^h \|_{L^2}$ &  $\frac{|\tilde \lambda^h_6-\lambda_6|}{\lambda_6}$ &  $|u_6 -\tilde u^h_6 |_{H^1}$ & $\|u_6 - \tilde u^h_6\|_{L^2}$ \\[0.1cm] \hline
	& 5	& 3.52e-7	& 4.95e-3	& 1.67e-4	& 3.05e-1	& 1.83e1	& 9.71e-1 \\[0.1cm]
	& 10	& 1.32e-9	& 5.75e-4	& 9.25e-6	& 3.21e-3	& 1.50	& 3.38e-2 \\[0.1cm]
3	& 20	& 5.09e-12	& 7.05e-5	& 5.60e-7	& 9.57e-6	& 1.12e-1	& 1.01e-3 \\[0.1cm]
	& 40	& 1.12e-13	& 8.77e-6	& 3.47e-8	& 3.45e-8	& 1.20e-2	& 4.94e-5 \\[0.1cm]
	&  $\rho_3$       & 8.04	& 3.05	& 4.08	& 7.76	& 3.54	& 4.79 \\[0.1cm] \hline

	& 5	& 6.90e-9	& 5.26e-4	& 1.80e-5	& 3.05e-1	& 1.87e1	& 9.91e-1 \\[0.1cm]
4	& 10	& 6.31e-12	& 2.91e-5	& 4.72e-7	& 7.59e-4	& 6.38e-1	& 1.45e-2 \\[0.1cm]
	& 20	& 5.75e-14	& 1.76e-6	& 1.41e-8	& 4.42e-7	& 1.91e-2	& 1.75e-4 \\[0.1cm]
	&  $\rho_4$    & 10.09	& 4.11	& 5.16	& 9.70	 & 4.97	& 6.23 \\[0.1cm] \hline

	& 5	& 1.13e-10	& 5.65e-5	& 1.97e-6	& 3.06e-1	& 1.89e1	& 1.00 \\[0.1cm]
5	& 10	& 1.15e-14	& 1.48e-6	& 2.43e-8	& 1.41e-4	& 2.63e-1	& 6.17e-3 \\[0.1cm]
	& 20	& 1.27e-14	& 4.42e-8	& 3.54e-10	& 1.72e-8	& 3.30e-3	& 3.06e-5 \\[0.1cm]
	&  $\rho_5$   & 13.26	& 5.16	& 6.22	& 12.04	& 6.24	& 7.50 \\[0.1cm] \hline
\end{tabular}
\end{table}

\begin{table}[ht]
\caption{Errors for the first and sixth eigenvalues in 2D and 3D when using IGA with optimally-blended quadratures and the boundary penalty technique.}
\label{tab:eig2d3d} 
\centering 
\begin{tabular}{| c | c || cc |cc | cc | cc |}
\hline
 & $p$ &  \multicolumn{2}{c|}{3} & \multicolumn{2}{c|}{4}  &  \multicolumn{2}{c|}{5} \\[0.1cm] \hline
$d$ & $N$ & $\frac{|\tilde \lambda^h_1-\lambda_1|}{\lambda_1}$ &  $\frac{|\tilde \lambda^h_6-\lambda_6|}{\lambda_6}$ & $\frac{|\tilde \lambda^h_1-\lambda_1|}{\lambda_1}$ &  $\frac{|\tilde \lambda^h_6-\lambda_6|}{\lambda_6}$ & $\frac{|\tilde \lambda^h_1-\lambda_1|}{\lambda_1}$ &  $\frac{|\tilde \lambda^h_6-\lambda_6|}{\lambda_6}$   \\[0.1cm] \hline
	& 3	& 2.28e-5	& 1.63e-4	& 1.32e-6	& 8.52e-3	& 6.55e-8	& 1.30e-7 \\[0.1cm]
	& 6	& 8.05e-8	& 6.10e-4	& 1.09e-9	& 9.03e-5	& 1.29e-11	& 1.15e-5 \\[0.1cm]
2	& 12	& 3.07e-10	& 1.94e-6	& 1.07e-12	& 6.06e-8	& 5.21e-14	& 1.61e-9 \\[0.1cm]
	& 24	& 1.31e-12	& 7.15e-9	& 1.35e-14	& 5.37e-11	& 9.45e-14	& 2.86e-13 \\[0.1cm]
	& $\rho_p $ & 8.02	& 8.19	& 10.12	& 10.54	& 7.95	& 12.8 \\[0.1cm] \hline
	& 2	& 6.78e-4	& 3.33e-1	& 1.00e-4	& 3.33e-1	& 1.28e-5	& 3.33e-1 \\[0.1cm]
	& 4	& 2.15e-6	& 6.03e-4	& 6.74e-8	& 8.92e-5	& 1.79e-9	& 1.14e-5 \\[0.1cm]
3	& 8	& 7.94e-9	& 1.92e-6	& 5.96e-11	& 5.99e-8	& 5.30e-13	& 1.59e-9 \\[0.1cm]
	& 16	& 3.06e-11	& 7.07e-9	& 3.00e-15	& 5.30e-11	& 4.66e-15	& 3.66e-13 \\[0.1cm]
	& $\rho_p $  & 8.13	& 8.48	& 11.5	& 10.82	& 12.26	& 13.19 \\[0.1cm] \hline
\end{tabular}
\end{table}

\begin{figure}[h!]
\centering
\includegraphics[height=5cm]{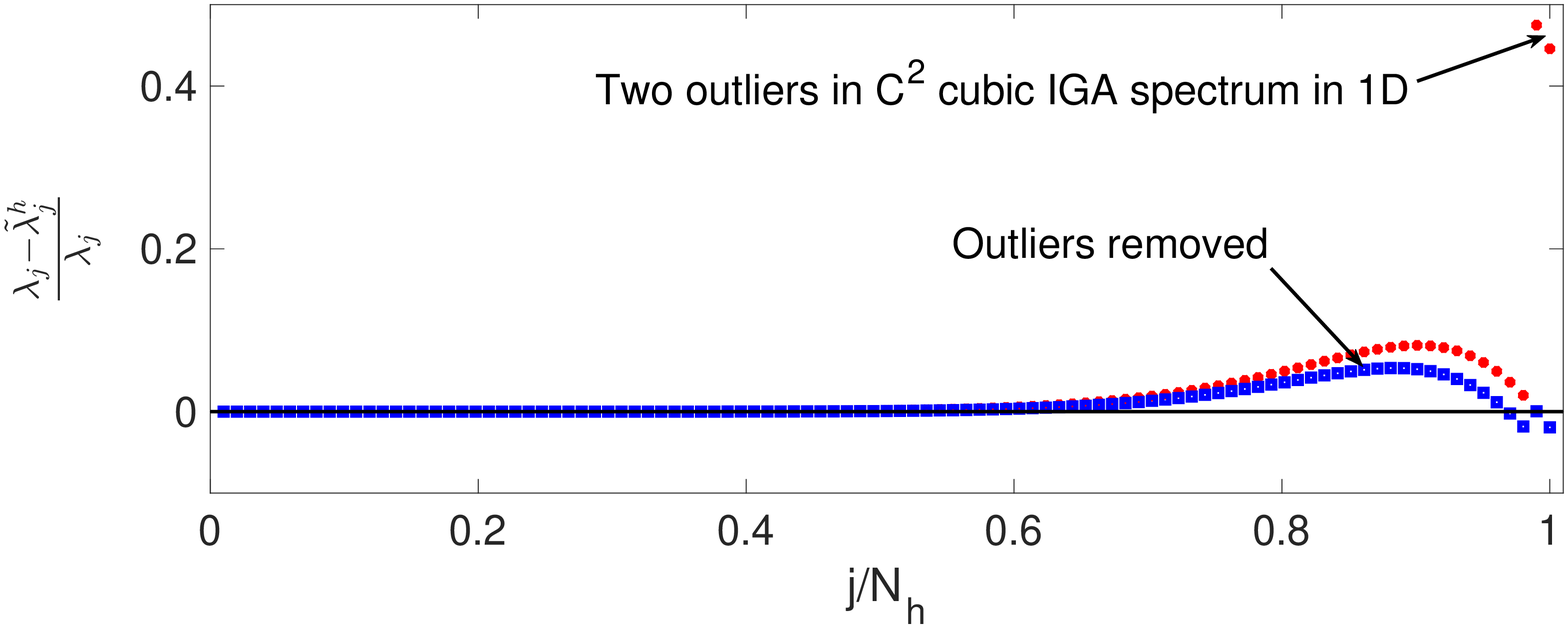} \\
\includegraphics[height=5cm]{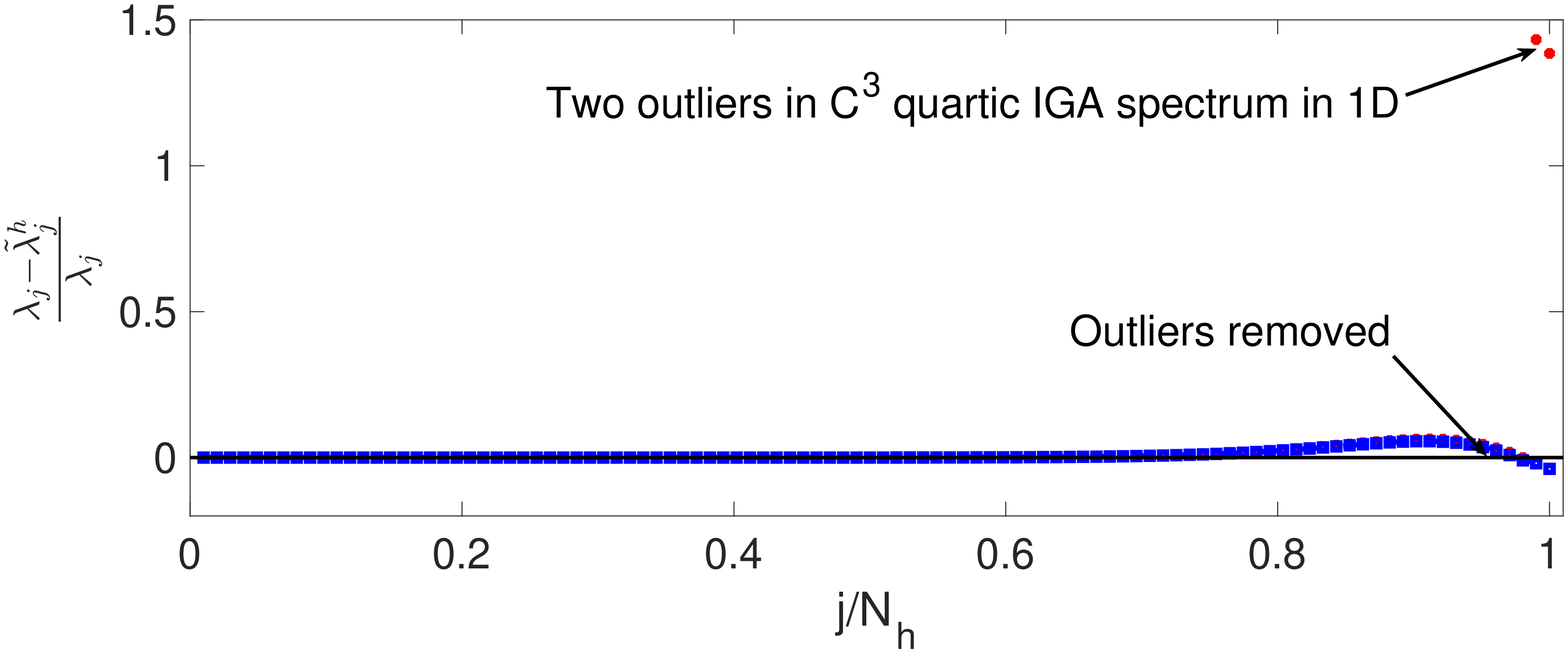} \\
\includegraphics[height=5cm]{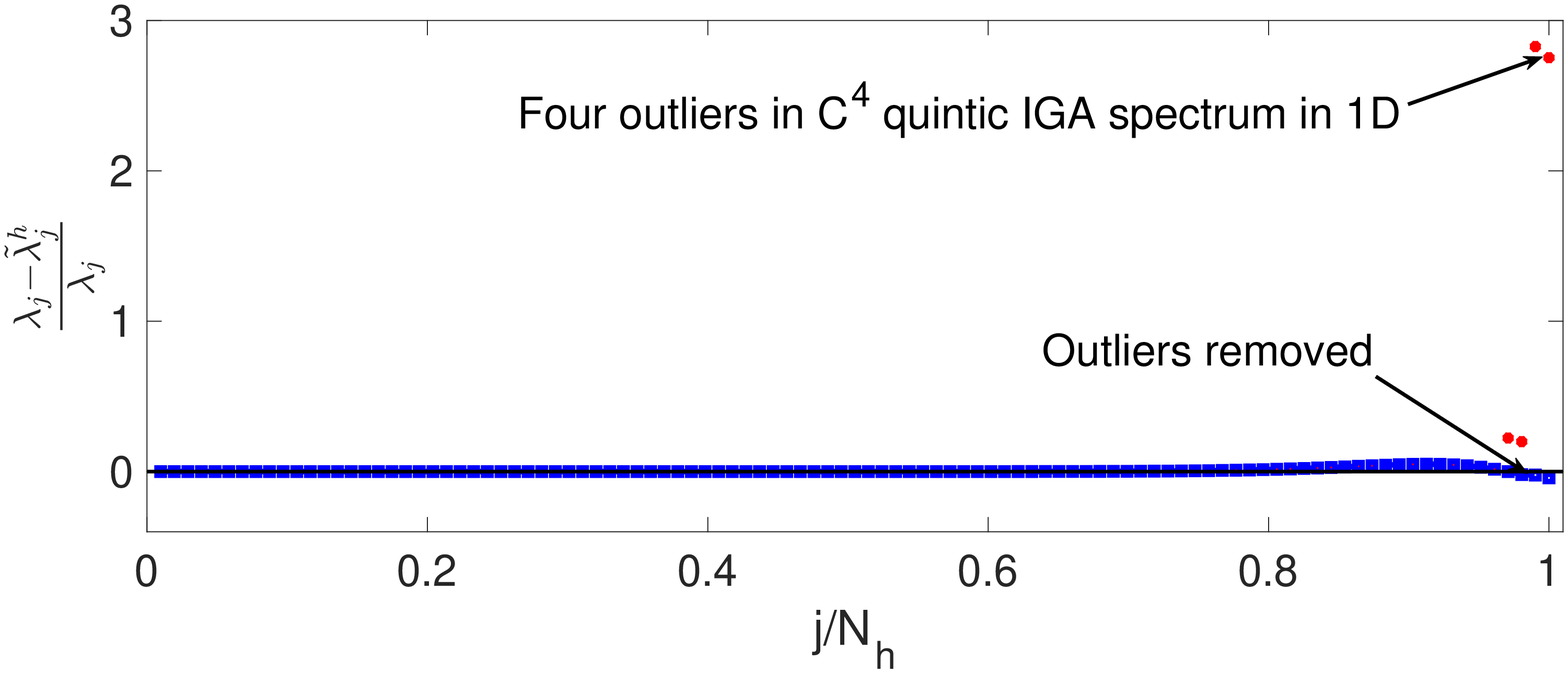} 
\vspace{-0.8cm}
\caption{Outliers in IGA spectra and their eliminations when using the IGA with optimally-blended rules and the boundary penalty technique. There are 100 elements in 1D with polynomial degrees $p\in\{3,4,5\}$.}
\label{fig:eig1d}
\end{figure}

\begin{figure}[h!]
\centering
\includegraphics[height=5cm]{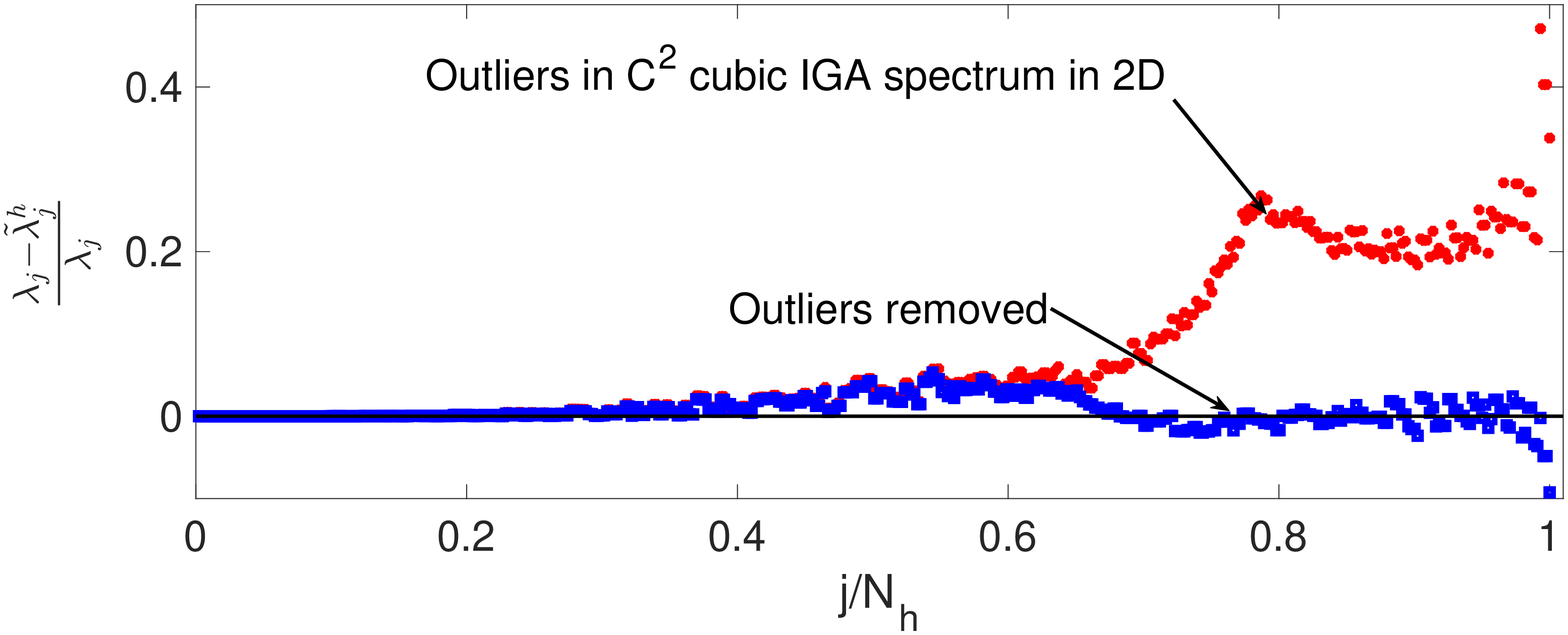} \\
\includegraphics[height=5cm]{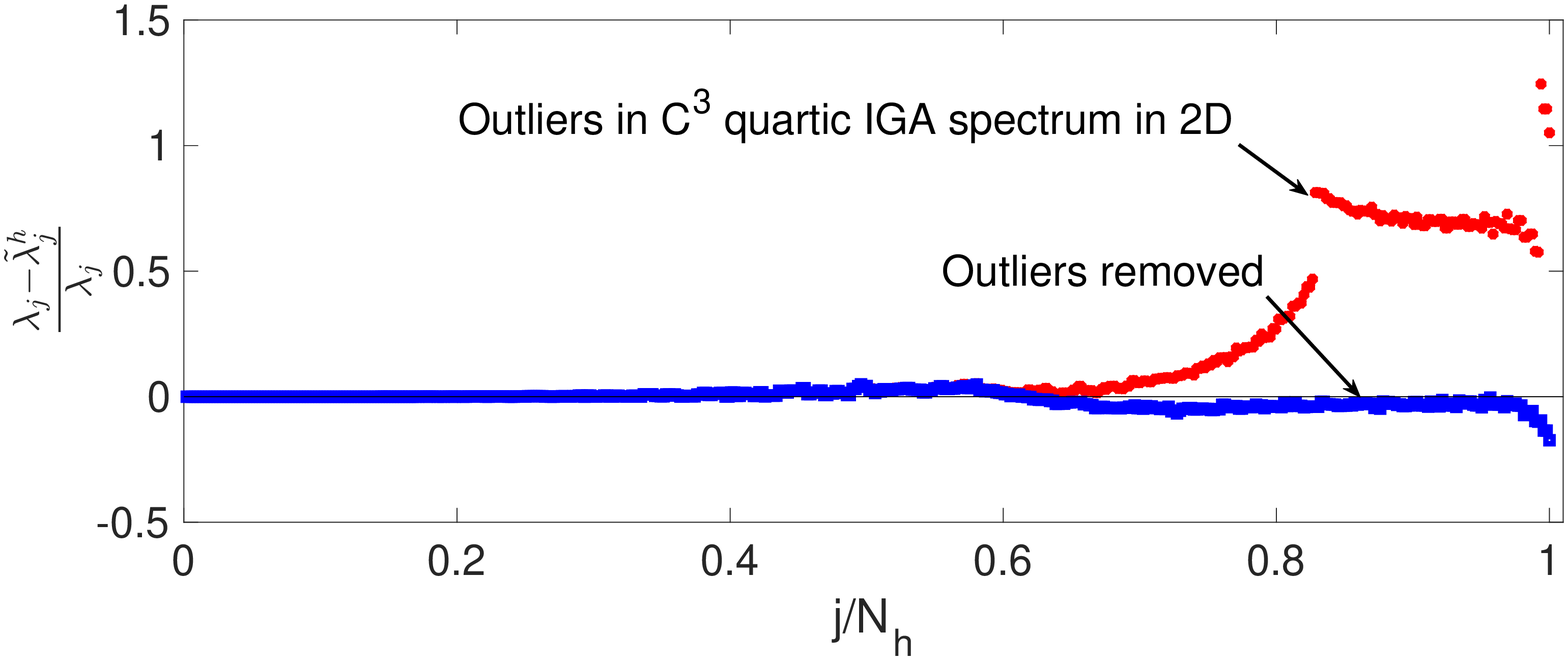} \\
\includegraphics[height=5cm]{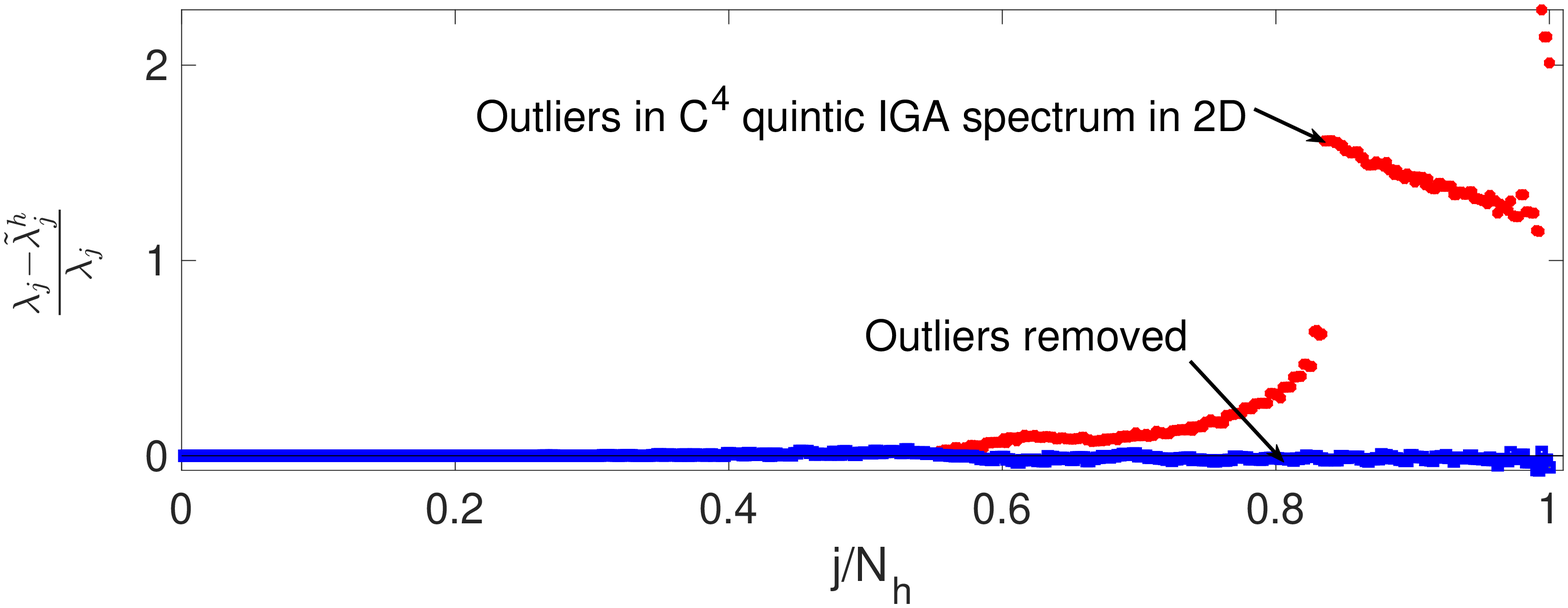} 
\vspace{-1cm}
\caption{Outliers in IGA spectra and their eliminations when using the IGA with optimally-blended rules and the boundary penalty technique. There are $20\times 20$ elements in 2D with polynomial degrees $p\in\{3,4,5\}$.}
\label{fig:eig2d}
\end{figure}

\begin{figure}[h!]
\centering
\includegraphics[height=5cm]{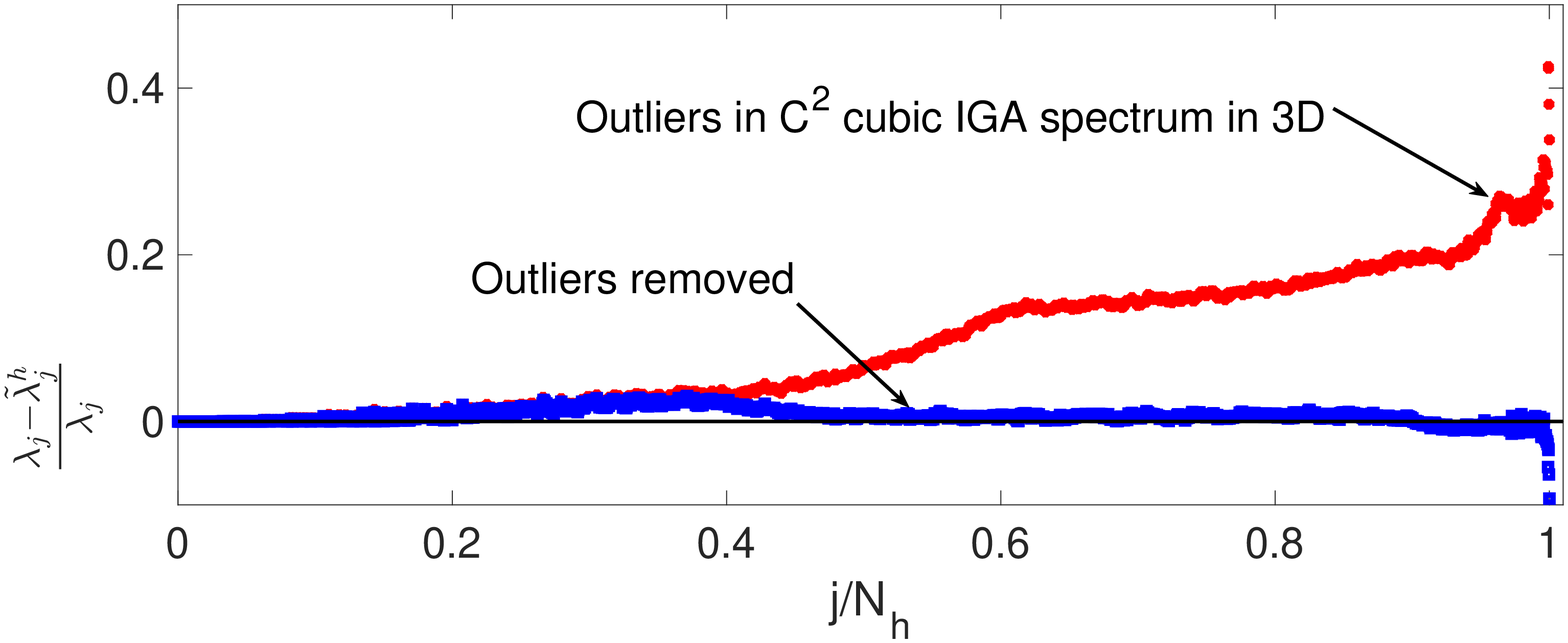} \\
\includegraphics[height=5cm]{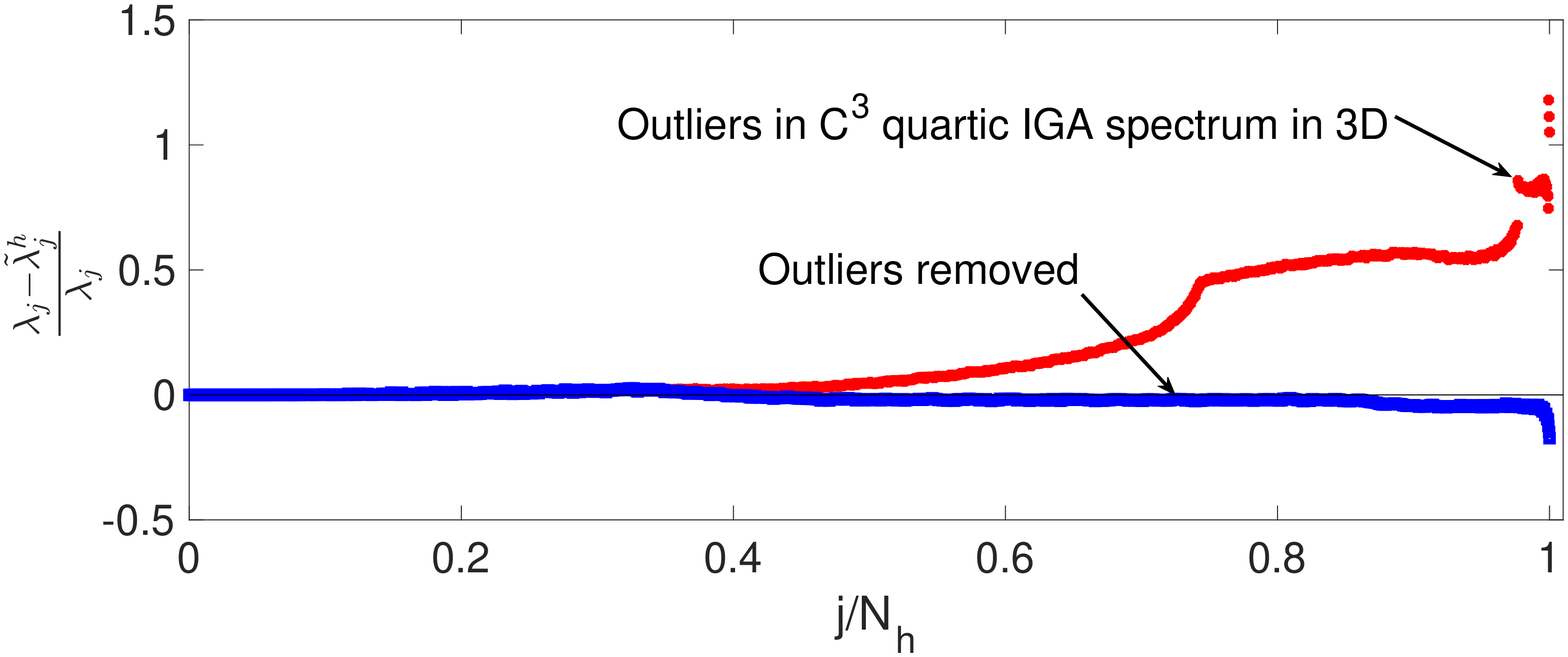} \\
\includegraphics[height=5cm]{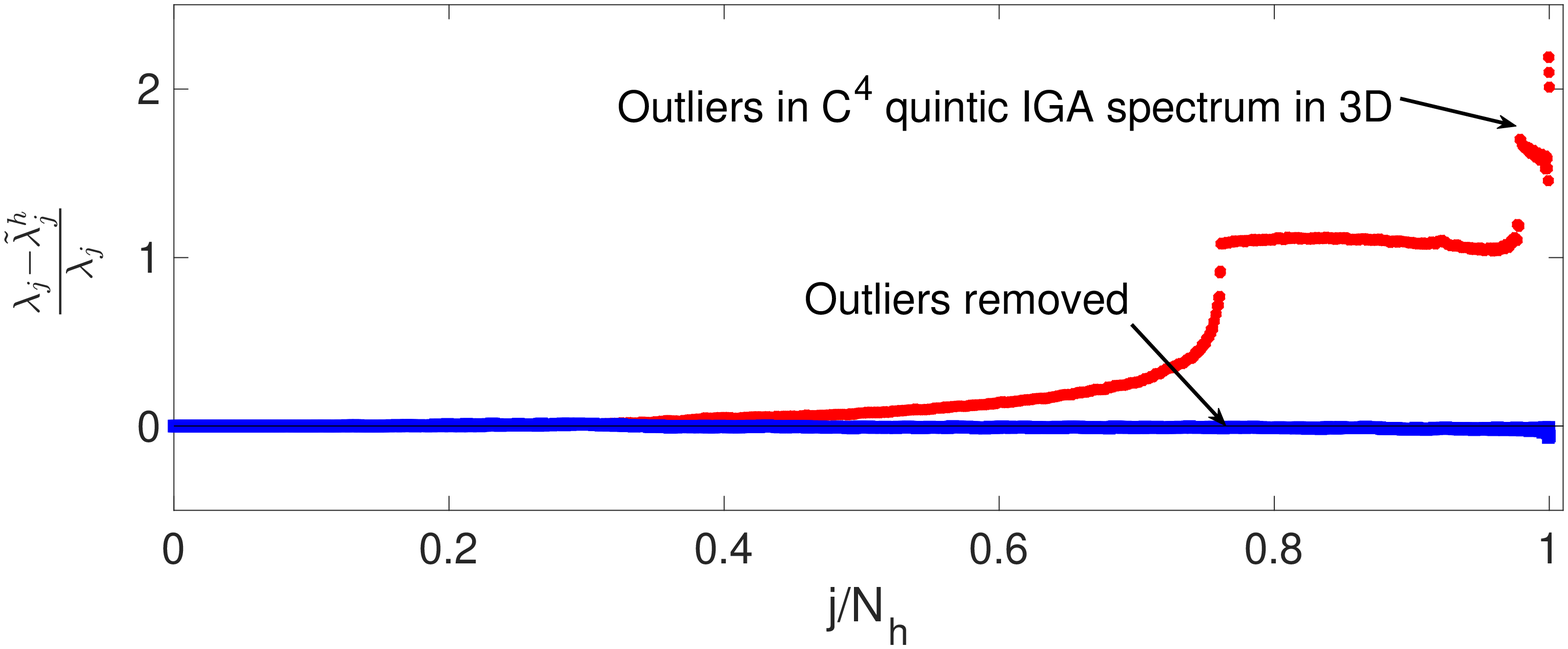} 
\vspace{-0.8cm}
\caption{Outliers in IGA spectra and their eliminations when using the IGA with optimally-blended rules and the boundary penalty technique. There are $20\times 20\times 20$ elements in 3D with polynomial degrees $p\in\{3,4,5\}$.}
\label{fig:eig3d}
\end{figure}

Figures \ref{fig:eig1d}--\ref{fig:eig3d} show the overall spectral errors when 
using the standard IGA 
and IGA with optimally-blended rules and the boundary penalty technique. 
The polynomial degrees are $p\in\{3,4,5\}$.
There are 100 elements in 1D, $20 \times 20$ elements in 2D, 
and $20\times 20\times 20$ elements in 3D, respectively. 
In all the scenarios, we observe that there are outliers in the IGA spectra. 
These outliers are eliminated by the proposed method. 
Moreover, we observe that these spectral errors are reduced significantly, especially in the high-frequency regions.

\subsection{Numerical study on condition numbers}
We study the condition numbers to show further advantages of the method. 
Since the stiffness and mass matrices are symmetric,
the condition numbers of the generalized matrix eigenvalue problems \eqref{eq:mevp} and 
\eqref{eq:mevp1d} are given by
\begin{equation}
\gamma := \frac{\lambda^h_{\max}}{\lambda^h_{\min}}, \qquad \tilde \gamma := \frac{\tilde \lambda^h_{\max}}{\tilde \lambda^h_{\min}},
\end{equation}
where $\lambda^h_{\max}, \tilde \lambda^h_{\max}$ are the largest eigenvalues and $\lambda^h_{\min}, \tilde \lambda^h_{\min}$ are the smallest eigenvalues of IGA and the proposed method, respectively. 
The condition number characterizes the stiffness of the system. 
We follow the recent work of soft-finite element method (SoftFEM) \cite{deng2020softfem} 
and define the \textit{condition number reduction ratio} of the method with respect to IGA as 
\begin{equation} \label{eq:srr}
\rho := \frac{\gamma}{\tilde\gamma} = \frac{\lambda^h_{\max} }{\tilde \lambda^h_{\max} } \cdot \frac{\tilde \lambda^h_{\min} }{\lambda^h_{\min}}.
\end{equation}
In general, one has $\lambda^h_{\min} \approx \tilde \lambda^h_{\min}$ for IGA and the proposed method with sufficient number of elements (in practice, these methods with only a few elements already lead to good approximations to the smallest eigenvalues). 
Thus, the condition number reduction ratio is mainly characterized by the ratio of the largest eigenvalues. 
Finally, we define the \textit{condition number reduction percentage} as
\begin{equation}
\varrho = 100 \frac{\gamma-\tilde \gamma}{\gamma}\,\% = 100(1-\rho^{-1}) \, \%.
\end{equation}

Table \ref{tab:cond} shows the smallest and largest eigenvalues, condition numbers and their reduction ratios 
and percentages for 1D, 2D, and 3D problems. 
We observe that the  condition numbers of the proposed method are significantly smaller. 
For higher-order isogeometric elements, there are more outliers and these outliers pollute a larger high-frequency region. Consequently, they lead to larger errors in the high-frequency region.  The proposed method reduces more errors for higher-order elements. This leads to small errors in the high-frequency region also for high-order elements. 
The condition number of the proposed method reduces by about 32\% for $C^2$ cubic, 60\% for $C^3$ quartic, and 75\% for $C^4$ quintic elements. This holds valid for both 2D and 3D problems.

\begin{table}[ht]
\caption{Minimal and maximal eigenvalues, condition numbers, reduction ratios and percentages when using IGA and IGA with optimally-blended quadratures and the boundary penalty technique. The polynomial degrees are $p\in\{3,4,5\}$. There are 100, $48\times 48$ elements,  and $16 \times 16 \times 16$ elements in 1D, 2D, and 3D, respectively.}
\label{tab:cond} 
\centering 
\begin{tabular}{|c| c | ccc | ccc c| cc |}
\hline
$d$ & $p$ & $\lambda_{\min}^h$ &  $\lambda_{\max}^h$ & $ \tilde \lambda_{\max}^h $ &  $\gamma$ &  $\tilde \gamma$ & $\rho$ & $\varrho$ \\[0.1cm] \hline
	& 3	& 9.87	& 1.46e5	& 9.87e4	& 1.47e4	& 1.00e4	& 1.47	& 32.17\% \\[0.1cm]
1	& 4	& 9.87	& 2.45e5	& 9.87e4	& 2.48e4	& 1.00e4	& 2.48	& 59.69\% \\[0.1cm]
	& 5	& 9.87	& 3.93e5	& 1.00e5	& 3.98e4	& 1.02e4	& 3.92	& 74.47\% \\[0.1cm] \hline
								
	& 3	& 1.97e1	& 6.71e4	& 4.55e4	& 3.40e3	& 2.30e3	& 1.47	& 32.17\% \\[0.1cm]
2	& 4	& 1.97e1	& 1.13e5	& 4.55e4	& 5.72e3	& 2.30e3	& 2.48	& 59.69\% \\[0.1cm]
	& 5	& 1.97e1	& 1.81e5	& 4.57e4	& 9.17e3	& 2.31e3	& 3.96	& 74.77\% \\[0.1cm] \hline
										
	& 3	& 2.96e1	& 1.12e4	& 7.58e3	& 3.78e2	& 2.56e2	& 1.48	& 32.23\% \\[0.1cm]
3	& 4	& 2.96e1	& 1.88e4	& 7.58e3	& 6.36e2	& 2.56e2	& 2.48	& 59.72\% \\[0.1cm]
	& 5	& 2.96e1	& 3.02e4	& 7.59e3	& 1.02e3	& 2.56e2	& 3.98	& 74.89\% \\[0.1cm] \hline
 \end{tabular}
\end{table}

\section{Concluding remarks} \label{sec:conclusion}
We improve the isogeometric spectral approximations by combining the two ideas: 
optimally-blended quadratures and a boundary penalty technique. 
As a result, we obtained a superconvergence of rate $h^{2p+2}$ for the eigenvalue errors 
and eliminated  the outliers in the spectra. 
The technique can be also used to improve the spectral approximations of the Neumann eigenvalue problems. 
These improvements lead to a better spatial discretization for the time-dependent partial differential equations, 
which in return, improve the overall performance of numerical methods. 
As future work, it would be interesting to study the method for higher-order differential operators 
and nonlinear application problems.

%
%
%
\bibliographystyle{splncs04}
\bibliography{igaref}
%
%
%
%
%
\end{document}